\documentclass{article}
\usepackage{amsfonts}
\usepackage{amsmath}
\usepackage[all]{xy}
\newtheorem{theorem}{Theorem}[section]
\newtheorem{definition}[theorem]{Definition}

\newtheorem{lemma}[theorem]{Lemma}
\newtheorem{algorithm}[theorem]{Algorithm}

\usepackage{mathrsfs}
\usepackage{amssymb}
\usepackage{amscd}
\date{}

\usepackage{amssymb}
\usepackage{extarrows}
\usepackage{cite}

\begin{document}

\title{Free Rota-Baxter systems and a Hopf algebra structure\footnote{Supported by the NNSF of China (no. 11571121),
the NSF of Guangdong Province (no. 2017A030313002) and the Science and Technology Program of Guangzhou (no. 201707010137).}}

\author{
Jianjun Qiu   \\
{\small \ School of Mathematics and Statistics, Lingnan Normal
University}\\
{\small Zhanjiang 524048, P. R. China}\\
{\small jianjunqiu@126.com}\\
\\
 Yuqun Chen\footnote {Corresponding author.}   \\
{\small \ School of Mathematical Sciences, South China Normal
University}\\
{\small Guangzhou 510631, P. R. China}\\
{\small yqchen@scnu.edu.cn}
}

\maketitle \noindent\textbf{Abstract:}
In this paper, we  give a linear basis of   a free    Rota-Baxter system on  a set   by using the  Gr\"{o}bner-Shirshov bases  method and then  we obtain a left counital  Hopf algebra structure on  a  free  Rota-Baxter system.

\ \


  \ \



\section{Introduction}

A triple $(A; R, S)$ consisting of an associative  unitary algebra $A$ over a field $k$ and two $k$-linear
operators $R, S : A \rightarrow A$ is called a Rota-Baxter system if, for any $a, b\in A$,
$$
R(a)R(b) = R(R(a)b + aS(b)), \ \ S(a)S(b) = S (R(a)b + aS(b)).
$$
Rota-Baxter system was   introduced by   Brzezi\'{n}ski in a recent paper \cite{Brz2016}, which  can be viewed as   an extension of the Rota-Baxter algebra of weight $0$.

An associative  unitary algebra  $A$ together with a $k$-linear
operator $ P:A\rightarrow A$ is called a   Rota-Baxter algebra of weight $\lambda$,  if
$$
P(x)P(y)  = P(  P(x)y)  +P( xP(y))+\lambda P( xy ),  \    x,y \in A.
$$
 Here $\lambda$ is a fixed element in the field $k$.

 Baxter \cite{Bax1960} firstly   studied   the    Rota-Baxter algebras.  Some combinatoric properties of Rota-Baxter  algebras were studied by
Rota  \cite{Rot1969} and  Cartier  \cite{Ca72}. The  constructions of   free   Rota-Baxter associative algebras  on both commutative and noncommutative cases  were given   by using  different methods, for example, \cite{  bcq, Ca72,   EG08a, EG08b,    lguo,    GK00a, GK00b, Rot1969}.

Gr\"{o}bner bases and Gr\"{o}bner-Shirshov bases  were invented
independently by   Shirshov for ideals of free
Lie algebras \cite{Sh} and implicitly free associative
algebras \cite{Sh}  (see also \cite{Be78,  Bo76}), by
Hironaka \cite{Hi64} for ideals of the power series algebras (both
formal and convergent), and by Buchberger \cite{Bu70} for ideals
of the polynomial algebras. Gr\"{o}bner bases and
Gr\"{o}bner-Shirshov bases theories have been proved to be very
useful in different branches of mathematics. See, for example, the books \cite{AL, BKu94,
BuW }, the papers \cite{Be78,  Bo76  }, and the
surveys \cite{BC14,   BK03}.  In fact,  Gr\"{o}bner-Shirshov bases theory   is a useful tool for constructing free objects for many  algebra varieties.   We will  apply  Gr\"{o}bner-Shirshov bases method to construct free Rota-Baxter system generated by a set.

The Hopf algebra originated from the study of topology and has widely applications on mathematics and physics. See for instance \cite{abe1980,Ck1998}.   Many classical  Hopf algebras are  build from free objects on various context, which include  the free associative algebra and the enveloping algebra of Lie algebras. Recently, there are some Hopf algebra structures on other  free objects, such as the dendriform algebras \cite{lo2004} and   Rota-Baxter algebras \cite{EFg2006, zgg2016, glt }. Inspired by the ideas of  the above papers, we will establish a left counital Hopf algebra structure on the  free  object of Rota-Baxter system.

The paper is organized as follows. In Section 2,    we give  a   Gr\"{o}bner-Shirshov basis  of  a  free    Rota-Baxter system  and then a  linear basis of  such    algebra  is  obtained  by  Composition-Diamond lemma for  associative  $\Omega$-algebras.
In Section 3, by using the construction of free Rota-Baxter system obtained in    Section 2,  we give a left counital  Hopf algebra structure on a free    Rota-Baxter system.

\section{Free    Rota-Baxter  systems}

\subsection{Gr\"{o}bner-Shirshov
bases   for  associative $\Omega$-algebras}

In this subsection, we review    Gr\"{o}bner-Shirshov
bases theory  for  associative unitary  $\Omega$-algebras. For more details, see for instance,     \cite{bcq, gsz2013}.

Let
$$
\Omega=\bigcup_{m=1}^{\infty}\Omega_{m},
$$
where $\Omega_{m}$ is a set of $m$-ary operators for any $m\geq 1$. For any set
 $Y$, let
$$
\Omega(Y)=\bigcup_{m=1}^{\infty}\left\{ \omega^{(m)}(y_1, y_2, \cdots, y_m)|y_i\in Y,  1 \leq i \leq m,  \omega^{(m)}\in \Omega_{m} \right\}.
$$
Let $X$ be a set. Define
$
\langle \Omega;  X \rangle_0=X^*$, where $X^*$ is the free monoid with the unit $1$ on the set $X$.    Assume that we have defined $\langle \Omega;  X \rangle_{n-1}$. Define
$$
\langle \Omega;  X \rangle_n= (X\cup \Omega(\langle \Omega;  X \rangle_{n-1}))^*.
$$
Then it is easy to see that
$
\langle \Omega;  X \rangle_{n }\subseteq \langle \Omega;  X \rangle_{n+1}
$
for any $n\geq 0$.
Set
$$
\langle \Omega;  X \rangle = \bigcup_{n=0}^{\infty}\langle \Omega;  X \rangle_n.
$$

If $u\in  \langle \Omega;  X \rangle$, then $u$ is said to be an $\Omega$-word  on the set $X$.
For any    $u \in X \cup \Omega( \langle \Omega;  X \rangle)$,   $u$  is called   prime. Therefore, for any $u\in \langle \Omega;  X \rangle$, $u$ can be uniquely expressed in the canonical form
$$
u=u_1u_2\cdots u_n,\ n\geq 0,
$$
where each $u_i$ is prime. The  breath of $u$, denoted by bre($u$), is defined to be  the number  $n$.     By the     depth of $u$, denoted by  dep($u$), we mean
$
{\rm dep}(u)=\min\{n|u\in \langle \Omega;  X \rangle_n \}.
$

An associative unitary $\Omega$-algebra  over a field $k$   is an associative unitary $k$-algebra $A$  with  a set of  multilinear operators $\Omega$, where $\Omega=\bigcup_{m=1}^{\infty}\Omega_{m}$
and each $\Omega_{m}$ is a set of $m$-ary multilinear operators on $A$.

Let $k\langle \Omega;  X\rangle$ be the linear space spanned by $\langle \Omega;  X\rangle$ over the field $k$. Then  $k\langle \Omega;  X\rangle$ is a  free   associative   unitary  $\Omega$-algebra on $X$.

Let    $\star\notin X$. By a
$\star$-$\Omega$-word on $X$ we mean any expression in $\langle \Omega; X\cup
\{\star\}\rangle$ with only one occurrence of $\star$.   Let  $\langle \Omega;  X \rangle^\star$ denote   the  set of all
  $\star$-$\Omega$-words on $X$. If  $\pi$ is a   $\star$-$\Omega$-word and $s\in k\langle
\Omega;  X\rangle$, then   we call
$
\pi|_{s}:=\pi|_{\star\mapsto s}
$
an   $s$-$\Omega$-word.

Now, we assume that $\langle \Omega;  X\rangle$ is equipped with a monomial
order  $>$. This means that $>$ is a well order  on
$\langle \Omega;  X\rangle$ such that for any $ v, w \in \langle \Omega;  X\rangle$ and
$\pi\in \langle \Omega;  X\rangle^\star$, if $w> v$, then  $\pi|_w> \pi|_v$.

For any $f\in k\langle \Omega;  X\rangle $, let
$\bar{f}$  be the leading   $\Omega$-word of $f$ with respect to the order $> $. If the coefficient
of $\bar{f}$ is $1$, then we call that $f$ is  monic. We also call  a  set $ \mathbb{S} \subseteq k\langle \Omega;  X\rangle $   monic if each  $s\in \mathbb{S}$ is monic.

\ \

Let $f, g\in k\langle \Omega;  X\rangle $ be   monic. Then we define   two
kinds of compositions.
\begin{enumerate}
\item[(I)]If   $w=\bar{f}a=b\bar{g}$ for some $a,b\in
\langle \Omega;  X\rangle$ such that bre($w$)$< $bre($\bar{f}$)$+$bre($\bar{g}$), then
we call $(f,g)_{w}=fa-bg$ the intersection composition of $f$
and $g$ with respect to the ambiguity $w$.
\item[(II)] If $w=\bar{f}=\pi|_{_{\bar{g}}}$ for some
$\pi \in \langle \Omega;  X\rangle^\star$, then we call $(f,g)_{w}=f-\pi|_{_{g}}$ the
inclusion  composition of  $f$ and $g$ with respect to the ambiguity $w$.
\end{enumerate}

If $\bar{f}=\pi|_{_{\bar{g}}}$, then the
  transformation $f\mapsto f-\alpha\pi|_{_{g}}$ is called the
Elimination of the Leading Word (ELW) of $f$ by $g$, where $g$ is monic and $\alpha$ is the coefficient of $\bar{f}$.

Let $\mathbb{S} \subseteq k\langle \Omega;  X\rangle$ be   monic. The composition
$(f,g)_w$ is called trivial modulo $(\mathbb{S},w)$ if
$$
(f,g)_w=\sum\alpha_i\pi_i|_{s_i},
$$
where each $\alpha_i\in k$,  $\pi_i\in \langle \Omega;  X\rangle^\star$, $s_i\in
\mathbb{S}$ and $\pi_i|_{\overline{s_i}}< w$. If this is the case, we write
$$
(f,g)_w\equiv 0 \ {\rm mod} (\mathbb{S},w).
$$

In general, for any   $p$ and $q$, $ p\equiv
q \ {\rm mod} (\mathbb{S},w) $ means that $ p-q=\sum\alpha_i\pi_i|_{s_i}, $ where
each $\alpha_i\in k$,  $\pi_i\in \langle \Omega;  X\rangle^\star$, $s_i\in \mathbb{S}$
and $\pi_i|_{\overline{s_i}}< w$.

A monic set $\mathbb{S}$ is called a Gr\"{o}bner-Shirshov basis  in  $k\langle
\Omega;  X\rangle$ if any composition $(f,g)_w$ of $f,g\in \mathbb{S}$  is
trivial modulo $(\mathbb{S},w)$.\\

In fact, the proof of  Composition-Diamond lemma  for   associative unitary $\Omega$-algebras is the same as the one for nonunitary $\Omega$-algebras in \cite{bcq}.  See also \cite{gsz2013}.

\begin{theorem}\label{cdla} {\em(\cite{bcq}, Composition-Diamond lemma  for   associative unitary $\Omega$-algebras)}\ \  Let $\mathbb{S} \subseteq k\langle \Omega;  X\rangle$ be   monic,   $> $ a monomial order  on $\langle \Omega;  X\rangle$  and  $Id(\mathbb{S})$ the $\Omega$-ideal of $k\langle \Omega;  X\rangle$ generated by $\mathbb{S}$.   Then the following
statements are equivalent:
 \begin{enumerate}
\item[(i)] The set $\mathbb{S} $ is a Gr\"{o}bner-Shirshov basis in $k\langle \Omega;  X\rangle$.
\item[(ii)] If  $ f\in Id(\mathbb{S})$, then $ \bar{f}=\pi|_{\overline{s}}$
for some $\pi \in \langle \Omega;  X\rangle^\star$ and $s\in \mathbb{S}$.
\item[(iii)]The set $Irr(\mathbb{S}) = \{ w\in \langle \Omega;  X\rangle |  w \neq
\pi|_{\overline{s}},
  \pi \in \langle \Omega;  X\rangle^\star,    s\in \mathbb{S}\}$
is a  $k$-linear basis of   $k\langle \Omega;  X|\mathbb{S}\rangle:=k\langle
\Omega;  X\rangle/Id(\mathbb{S})$.
\end{enumerate}
\end{theorem}

\subsection{Free  Rota-Baxter systems}
In this subsection,  we give  a   Gr\"{o}bner-Shirshov basis  of  a  free    Rota-Baxter system on a set  and then a  linear basis of  such  an  algebra  is  obtained  by the Composition-Diamond lemma for  associative   unitary $\Omega$-algebras.\\

Let $X$   be a well-ordered set and $\Omega=\{R, S\}$, where both  $R$ and  $S$  are    1-array operators. Assume that $R>S$.   If $u\in \langle \{R, S\},  X\rangle$, then we define ${\rm deg}(u)$ to be the number of   all   occurrences of all $x\in X$ and  $\omega \in \Omega$.  For example, if $w=x_1x_2R(R(x_1)S(x_2))$, where $x_1, x_2\in X$, then   ${\rm deg}(u)=7$.

If
$
u=u_1u_2\cdots u_n\in \langle \{R, S\};  X\rangle, n\geq 1,
$
where each $u_i$ is prime, then we let
$$
{\rm wt}(u)=(\mbox{deg}(u), \mbox{bre}(u), u_1,u_2,\cdots, u_n).
$$
  Define the  Deg-lex order  $ >_{_{{\rm Dl}}} $ on $\langle \{R, S\};  X\rangle$ as follows. For any $u=u_1u_2\cdots u_n$, $v=v_1v_2\cdots v_m\in \langle \{R, S\};  X\rangle$,   where    $u_i,  v_j$  are  prime, define
$$
u>_{_{{\rm Dl}}}v \ \mbox{if}\ {\rm wt}(u)>{\rm wt}(v)\  \mbox{lexicographically},
$$
where $u_i> v_i$ if  deg($u_i$)$>$deg($v_i$) or  deg($u_i$)=deg($v_i$) such that  one of the following conditions  holds:

(a)  $u_i, v_i\in X$ and $u_i>v_i$;

(b)  $u_i=R(1)$ or $u_i=S(1)$  and   $v_i\in X$;

(c)  $u_i=\omega (u_i'), v_i=\theta(v_i')$, $\omega, \theta\in \{R, S\}$ and
$$
 (\omega , u_i' )>(\theta, v_i')\  \mbox{lexicographically}.
$$

It is easy to see that $>_{_{{\rm Dl}}}$ is a monomial order on $\langle \{R, S\};  X\rangle$.

Let $k\langle \{R, S\}, X\rangle$ be the free associative  unitary  $\{R, S\}$-algebra generated by the set $X$.

\begin{theorem}\label{thgsb} With the order $>_{_{{\rm Dl}}}$ on $\langle \{R, S\};  X\rangle$,
the set
$$
\mathbb{S}=\left\{
\left.
 \begin{array}{ll}
 R(u)R(v)- R(R(u)v) - R(uS(v))\\
S(u)S(v)-\  S(R(u)v)- S(uS(v)) \\
\end{array}
\right|u, v\in \langle\{R, S\}, X\rangle \right\}
$$
is a Gr\"{o}bner-Shirshov basis in $k\langle \{R, S\}, X\rangle$.
\end{theorem}
{\bf Proof.} Let
$$
f(u,v)=R(u)R(v)- R(R(u)v) - R(uS(v)),
$$
$$
h(u, v)=S(u)S(v)- S(R(u)v)- S(uS(v)),
$$
where $u, \ v \in \langle\{R, S\}, X\rangle$.
  All the  possible compositions of
  $ \{R, S\}$-polynomials in $\mathbb{S}$ are listed as below:
$$
(f(u,v), f(v, w))_{w_1},  \ \  w_1=R(u)R(v)R(w),
$$
$$
(f(\pi|_{_{R(u)R(v)}}, w), f(u, v))_{w_2},\ \  w_2=R(\pi|_{_{R(u)R(v)}})R(w),
$$
$$
  f(u, \pi|_{_{R(v)R(w)}}), f(v, w))_{w_3},\ \  w_3=R(u)R(\pi|_{_{R(v)R(w)}}),
$$
$$
(f(\pi|_{_{S(u)S(v)}}, w), h(u, v))_{w_4},\ \  w_4=R(\pi|_{_{S(u)S(v)}})R(w),
$$
$$
 (f(u, \pi|_{_{S(v)S(w)}}), h(v, w))_{w_5},\ \  w_5=R(u)R(\pi|_{_{S(v)S(w)}}),
$$
$$
(h(u,v), h(v, w))_{w_6},  \ \ \  w_6=S(u)S(v)S(w),
$$
$$
(h(\pi|_{_{S(u)S(v)}}, w), h(u, v))_{w_7},\ \  w_7=S(\pi|_{_{S(u)S(v)}})S(w),
$$
$$
  (h(u, \pi|_{_{S(v)S(w)}}), h(v, w))_{w_8},\ \  w_8=S(u)S(\pi|_{_{S(v)S(w)}}),
$$
$$
(h(\pi|_{_{R(u)R(v)}}, w), f(u, v))_{w_9},\ \  w_9=S(\pi|_{_{R(u)R(v)}})S(w),
$$
$$
 (h(u, \pi|_{_{R(v)R(w)}}), f(v, w))_{w_{_{10}}},\ \  w_{_{10}}=S(u)S(\pi|_{_{R(v)R(w)}}).
$$
We check that all the compositions are trivial. Here, we just check one as an example.
\begin{eqnarray*}
&&(f(u,v), f(v, w))_{w_1}\\
&=& f(u, v)R(w)-R(u)f(v, w) \\
&=&  - (R(R(u)v)+R(uS(v)))R(w)+R(u)(R(R(v)w)+R(vS(w))) \\
&\equiv&  -R(R(R(u)v)w)-R(R(u)vS(w))-R(R(uS(v))w)-R(uS(v)S(w))\\
& &  +R(R(u)R(v)w)+R(uS(R(v)w))+R(R(u)vS(w))+R(uS(vS(w)))\\
&\equiv&  -R(R(R(u)v)w)-R(R(uS(v))w) -R(uS(R(v)w))-R(uS(vS(w)))\\
& &  +R(uS(R(v)w))+R(uS(vS(w)))+R(R(R(u)v)w)+R(R(uS(v))w)\\
&\equiv&\  0 \ {\rm mod} (\mathbb{S}, w_1).
\end{eqnarray*}
\hfill $ \square$\\

Note  that
$
RS(X)= k\langle \{R, S\}, X|\mathbb{S} \rangle
$
is a free  Rota-Baxter system on the set $X$, where
$$
\mathbb{S}=\left\{
\left.
\begin{array}{ll}
R(u)R(v)- R(R(u)v) - R(uS(v)) \\
S(u)S(v)-\  S(R(u)v)- S(uS(v)) \\
\end{array}
\right|u, v\in \langle\{R, S\}, X\rangle \right\}.
$$

Let
$$
\Phi_{\infty}(X) =\left\{w\in \langle \{R, S\}; X\rangle
\left|
\begin{array}{ll}
w\neq \pi|_{_{Q(u)Q(v)}},  \pi\in \langle \{R, S\}, X\rangle^{\star},   \\
  Q\in \{R, S\}, \ u,\   v\in \langle\{R, S\}; X\rangle\\
\end{array}
\right. \right\}.
$$
The elements of $\Phi_{\infty}(X)$ are called  Rota-Baxter system words.\\

\begin{theorem}
The set   $Irr(\mathbb{S})=\Phi_{\infty}(X)$ is a  linear basis  of the free   Rota-Baxter system  $RS(X)$.
\end{theorem}
{\bf Proof.} By   Theorems \ref{cdla} and \ref{thgsb}, we can obtain the result.\hfill $ \square$\\

By using ELWs, we have the following algorithm, which is    an  algorithm to
compute the product of two  Rota-Baxter system words in the free Rota-Baxter system  $RS(X)$.
 \\

 \noindent
\begin{algorithm}
  Let $w,v\in  \Phi_{\infty}(X)$. We define $w\diamond v$ by
induction on $n={\rm dep}(w)+{\rm dep}(v)$.
\begin{enumerate}
\item[(a)]  If $n=0$, then $w,v \in X^*$ and $w\diamond
v=wv$;
\item[(b)] If $n\geq 1$, there are two cases to consider:

(i) If  ${\rm bre}(w)={\rm bre}(v)=1$,   then
 \newline
\begin{equation*}
 w\diamond v=\left\{
\begin{array}{l@{\quad\quad}l}
Q(R(\widetilde{w})\diamond \widetilde{v} + \widetilde{w}\diamond S(\widetilde{v})),  &
{if} \  w=Q(\widetilde{w}),   v=Q(\widetilde{v}), Q\in \{R, S\},\\
wv, &   {  otherwise.}
\end{array}
\right.
\end{equation*}

(ii) If ${\rm bre}(w)>1$ or ${\rm bre}(v)>1$ and assume that $ w=w_1w_2\cdots w_t$  and
 $ v=v_1v_2\cdots v_l, $ where $u_i$ and $v_j$ are
prime, then
$$
 w\diamond v =w_1w_2\cdots w_{t-1} (w_t\diamond v_1)v_2\cdots v_l.
$$
\end{enumerate}

\end{algorithm}

\section{A  left counital Hopf algebra structure on free Rota-Baxter system}

In this section, similar to the   Hopf algebra  structure on free Rota-Baxter algebra given by Gao, Guo and Zhang\cite{glt},  we establish a left counital  Hopf algebra  structure on the free Rota-Baxter system $RS(X)$.
\subsection{A left counital  bialgebra structure}

In this subsection, we give a left counital  bialgebra structure on the free Rota-Baxter system $RS(X)$.

\begin{definition}{\em (\cite{zg2017})}
(a)\ \   A left counital coalgebra is a triple $(H, \Delta, \varepsilon)$,   where the comultiplication $\Delta: H\rightarrow H\otimes H$ satisfies the   coassociativity and the left counit  $\varepsilon: H\mapsto k$ satisfies the left
counicity: $
(\varepsilon\otimes{\rm id})\Delta =\beta_l,
$
where
$
\beta_l(w): RS(X)\rightarrow k\otimes RS(X), w\mapsto 1_k\otimes w.
$

(b)\ \  A left counital bialgebra is a quintuple  $(H, \mu, u, \Delta, \varepsilon)$, where $(H, \mu, u)$ is a $k$-algebra and  $(H, \Delta, \varepsilon)$ is a left counital coalgebra  such that $\Delta: H\rightarrow H \otimes H $ and $\varepsilon: H\rightarrow k$ are algebra homomorphisms.
\end{definition}

Define $\Delta:  RS(X) \rightarrow RS(X)\otimes RS(X)$, where for any $w\in \Phi_{\infty}(X)$, $\Delta(w)$ is defined by   induction  on   ${\rm dep}(w)$:

If ${\rm dep}(w)=0$,   then we define

\begin{equation*}
\Delta(w)=\left\{
\begin{array}{l@{\quad\quad}l}
1\otimes 1,  &
\mbox{if} \  w=1,\\
1\otimes x +x\otimes 1,  &\mbox{if}\  w=x\in X,
\end{array}
\right.
\end{equation*}
and
$$
\Delta(w)=\Delta(x_1)\diamond\Delta(x_2)\diamond \cdots \diamond \Delta(x_m),
$$
where $w=x_1x_2\cdots x_m, \ m\geq 2$ with each  $x_i\in X$.

Assume that $\Delta(w)$ has been defined for any $w\in \Phi_{\infty}(X)$ with ${\rm dep}(w)\leq n $.  Let $w\in \Phi_{\infty}(X)$ with  ${\rm dep}(w)=n+1$.  If ${\rm bre}(w)=1$, define

\begin{equation*}
\Delta(w)=\left\{
\begin{array}{l@{\quad\quad}l}
R(\widetilde{w})\otimes 1+({\rm id}\otimes R)\Delta(\widetilde{w}), &
\mbox{if} \  w=R(\widetilde{w}),\\
R(\widetilde{w})\otimes 1+({\rm id}\otimes S)\Delta(\widetilde{w}),  &\mbox{if}\  w=S(\widetilde{w}),
\end{array}
\right.
\end{equation*}
where $\Delta(\widetilde{w})$ is defined by induction hypothesis.  If ${\rm bre}(w)>1$, say $w=w_1w_2\cdots w_m$, where each $w_i$ is prime,  define
$$
\Delta(w)=\Delta(w_1)\diamond \Delta(w_2)\diamond \cdots \diamond\Delta(w_m).
$$
\ \

Define  $
u: k\rightarrow RS(X), 1_k\mapsto 1
$ and  $\varepsilon: RS(X) \rightarrow k$ by
\begin{equation*}
\varepsilon(w)=\left\{
\begin{array}{l@{\quad\quad}l}
1_k,  &
\mbox{if}  \   w=1, \\
0,  & \mbox{otherwise},
\end{array}
\right.
\end{equation*}
where $1_k$ is the unit of the field $k$.

It is easy to see that $(RS(X), \diamond, u)$ is a unitary $k$-algebra.

\begin{lemma}\label{lem3.2}
For any  $w, v \in RS(X)$, we have
$$
\Delta(w\diamond v)=\Delta(w)\diamond\Delta(v).
$$
\end{lemma}
{\bf Proof.} It is sufficient to prove that  $
\Delta(w\diamond v)=\Delta(w)\diamond\Delta(v)
$
for any $w, v\in \Phi_{\infty}(X)$. Induction on $(m, n)$, where $m={\rm dep}(w)+{\rm dep}(v)$ and  $n={\rm bre}(w)+{\rm bre}(v)$.

If $(m, n)=(0, 0)$, i.e. $w=v=1$, then $\Delta(w\diamond v)= 1\otimes1=\Delta(w)\diamond \Delta(v)$.

If $(m, n)=(0,1)$, then $w\in X,\ v=1$ or $w=1,\ v\in X$. Thus $\Delta(w\diamond v)=\Delta(w)\diamond\Delta(v)$.

Assume that we have proved   $\Delta(w\diamond v)=\Delta(w)\diamond\Delta(v)$  for all $u, v\in \Phi_{\infty}(X)$ with $(m, n)<(p, q)$.

Let $u, v\in \Phi_{\infty}(X)$  with $({\rm dep}(w)+{\rm dep}(v), {\rm bre}(w)+{\rm bre}(v))=(p, q)$ . There are two cases to consider.

 Case 1.   If   $q=1$, i.e.  $w =1$ or $v=1$, then by definition, the result is true.

 Case 2.   If   $q=2$,   then we have two cases.

(i) If $(w, v)\notin \bigcup _{Q\in \{R, \ S \}}Q(\Phi_{\infty}(X))\times Q(\Phi_{\infty}(X)) $, then $w\diamond v=wv\in \Phi_{\infty}(X)$ and the result is true by the definition of $\Delta$.

(ii) If $(w, v)\in \bigcup _{Q\in \{R, \ S \}}Q(\Phi_{\infty}(X))\times Q(\Phi_{\infty}(X)) $,  then   $w=R(w')$ and $v=R(v')$ or $w=S(w')$ and $v=S(v')$. Using the Sweedler notation, we can write
$$
\Delta(w')=\sum_{(w')} w'_{_{(1)}}\otimes w'_{_{(2)}}, \ \ \ \Delta(v')=\sum_{(v')} v'_{_{(1)}}\otimes v'_{_{(2)}}.
$$
(a) If $w=R(w')$ and $v=R(v')$, then
\begin{eqnarray*}
&&\Delta(R(w')\diamond R(v'))\\
&=&\Delta(R(R(w')\diamond v'+w'\diamond S(v')))\\
&=& R(R(w')\diamond v'+w'\diamond S(v'))\otimes 1+({\rm id}\otimes R)\Delta  (R(w')\diamond v'+w'\diamond S(v'))\\
&=& (w\diamond v)\otimes 1+({\rm id}\otimes R)\Delta  (R(w')\diamond v'+w'\diamond S(v'))\\
&=& (w\diamond v)\otimes 1+({\rm id}\otimes R)( \Delta(R(w'))\diamond \Delta (v')+ \Delta (w')\diamond \Delta (S(v')))\\
&=& (w\diamond v)\otimes 1+({\rm id}\otimes R)(  (R(w')\otimes 1+({\rm id}\otimes R)\Delta(w'))\diamond \Delta (v')\\
 && +\ \Delta (w')\diamond (R(v') \otimes 1 + ({\rm id}\otimes S)\Delta(w') ))\\
&=& (w\diamond v)\otimes 1+({\rm id}\otimes R)(\sum_{(v')}( R(w')\diamond  v'_{_{(1)}})  \otimes v'_{_{(2)}}
+ \sum_{(w')}(w'_{_{(1)}}\diamond R(v')) \otimes w'_{_{(2)}} \\
 && +\ \sum_{(w'),(v')}(w'_{_{(1)}}\diamond v'_{_{(1)}})\otimes (R(w'_{_{(2)}})\diamond v'_{_{(2)}}+w'_{_{(2)}}\diamond S(v'_{_{(2)}}) ))\\
&=& (w\diamond v)\otimes 1+  \sum_{(v')}( R(w')\diamond  v'_{_{(1)}})  \otimes R(v'_{_{(2)}}) +  \sum_{(w')}(w'_{_{(1)}}\diamond R(v')) \otimes R(w'_{_{(2)}}) \\
 && +\ \sum_{(w'),(v')}(w'_{_{(1)}}\diamond v'_{_{(1)}})\otimes (R(w'_{_{(2)}})\diamond R(v'_{_{(2)}}))\\
&=&(R(w')\otimes 1+\sum_{(w')}w'_{_{(1)}} \otimes R(w'_{_{(2)}}) )\diamond (R(v')\otimes 1+\sum_{(v')}v'_{_{(1)}} \otimes R(v'_{_{(2)}}))
\end{eqnarray*}
\begin{eqnarray*}
&=&(R(w')\otimes 1+({\rm id }\otimes R) (  \sum_{(w')} w'_{_{(1)}} \otimes   w'_{_{(2)}} ) )\diamond (R(v')\otimes 1+({\rm id }\otimes R) ( \sum_{(v')}v'_{_{(1)}} \otimes  v'_{_{(2)}}) )\\
&=&(R(w')\otimes 1+({\rm id }\otimes R) \Delta(w'))\diamond (R(v')\otimes 1+({\rm id }\otimes R)\Delta(v'))\\
&=& \Delta(w)\diamond \Delta(v).
\end{eqnarray*}
(b) If $w=S(w')$ and $v=S(v')$, then
\begin{eqnarray*}
&&\Delta(S(w')\diamond S(v'))\\
&=&\Delta(S(R(w')\diamond v'+w'\diamond S(v')))\\
&=& R(R(w')\diamond v'+w'\diamond S(v'))\otimes 1+({\rm id}\otimes S)\Delta (R(w')\diamond v'+w'\diamond S(v'))\\
&=& (R(w')\diamond R( v'))\otimes 1+({\rm id}\otimes S)( \Delta(R(w'))\diamond \Delta (v')+ \Delta (w')\diamond \Delta (S(v')))\\
&=& (R(w')\diamond R( v'))\otimes 1+({\rm id}\otimes S)(  (R(w')\otimes 1+({\rm id}\otimes R)\Delta(w'))\diamond \Delta (v')\\
 && +\  \Delta(w')\diamond (R(v') \otimes 1 + ({\rm id}\otimes S)\Delta(w') ))\\
&=& (R(w')\diamond R( v'))\otimes 1+({\rm id}\otimes S)(\sum_{(v')}( R(w')\diamond  v'_{_{(1)}})\otimes v'_{_{(2)}}
+ \sum_{(w')}(w'_{_{(1)}}\diamond R(v')) \otimes w'_{_{(2)}}\\ && +\ \sum_{(w'),(v')}(w'_{_{(1)}}\diamond v'_{_{(1)}})\otimes (R(w'_{_{(2)}})\diamond v'_{_{(2)}}+w'_{_{(2)}}\diamond S(v'_{_{(2)}}) ))\\
&=& (R(w')\diamond R( v'))\otimes 1+  \sum_{(v')}( R(w')\diamond  v'_{_{(1)}})  \otimes S(v'_{_{(2)}}) +  \sum_{(w')}(w'\diamond R(v')) \otimes S(w'_{_{(2)}}) \\
 && +\ \sum_{(w'),(v')}(w'_{_{(1)}}\diamond v'_{_{(1)}})\otimes (S(w'_{_{(2)}})\diamond S(v'_{_{(2)}}))\\
&=&\left(R(w')\otimes 1+({\rm id }\otimes S) \Delta(w') \right)\diamond \left(R(v')\otimes 1+({\rm id }\otimes S)\Delta(v')  \right)\\
&=& \Delta(w)\diamond \Delta(v).
\end{eqnarray*}

Case 3. If $q>2$, say $w=w_1w_2\cdots w_t$ and $v=v_1v_2\cdots v_l$, then we have two cases to consider.

(i) If $(w_t, v_1)\notin \bigcup _{Q\in \{R, \ S \}}Q(\Phi_{\infty}(X))\times Q(\Phi_{\infty}(X)) $, then $w\diamond v=wv\in \Phi_{\infty}(X)$ and the result is true by the definition of $\Delta$.

(ii)  If $(w_t, v_1)\in \bigcup _{Q\in \{R, \ S \}}Q(\Phi_{\infty}(X))\times Q(\Phi_{\infty}(X)) $, then by Case 2 or by induction, we have
$\Delta(w_t\diamond v_1)=\Delta(w_t)\diamond \bigtriangleup(v_1)$. By the definition of $\diamond$, we have  $w_t \diamond  v_1=\sum u_i$, where  ${\rm bre}(u_i)=1$. Therefore
\begin{eqnarray*}
\Delta (w\diamond  v)  &=& \Delta (w_1\diamond w_2\diamond \cdots \diamond (w_t \diamond  v_1)\diamond v_2 \diamond \cdots \diamond v_l)\\
&=& \sum  \Delta ( w_1\diamond w_2\diamond \cdots \diamond u_i\diamond v_2 \diamond \cdots \diamond v_l)\\
&=& \sum  \Delta ( w_1) \diamond \Delta(w_2)\diamond \cdots \diamond \Delta (u_i)\diamond \Delta(v_2) \diamond \cdots \diamond\Delta (v_l)\\
&=& \Delta ( w_1) \diamond \Delta(w_2)\diamond \cdots \diamond  (\sum  \Delta (u_i))\diamond \Delta(v_2) \diamond \cdots \diamond\Delta (v_l)\\
&=&\Delta ( w_1) \diamond \Delta(w_2)\diamond \cdots \diamond  \Delta(w_t)\diamond \Delta(v_1)\diamond \Delta(v_2) \diamond \cdots \diamond\Delta (v_l)\\
&=&\Delta (w)\diamond  \Delta(v).
\end{eqnarray*}
 This completes the proof by induction. \hfill $ \square$\\

\begin{lemma}\label{lem3.3}
For any  $w, v \in RS(X)$, we have
$
\varepsilon(w\diamond v)=\varepsilon(w)\diamond\varepsilon(v).
$
\end{lemma}
{\bf Proof.}  It is  easy to prove the result.\hfill $ \square$

\ \

Dr. Xing Gao has kindly pointed out that   $\varepsilon$ is not a right counit in   Lemma \ref{lem3.3}.

\begin{lemma}\label{lem3.4}
The triple $(RS(X),  \Delta, \varepsilon)$ is a  left couital  bialgebra.
\end{lemma}
{\bf Proof.} It is sufficient to  prove the coassociativity and the left counicity for $w\in \Phi_{\infty}(X)$.

We check the coassociativity by induction on $(m,n)$ for $w\in \Phi_{\infty}(X)$,    where $m={\rm dep}(w)$ and $n={\rm bre}(w)$.

 If $(m, n)=(0,0)$,  then $w=1$ and  $(\Delta\otimes {\rm id})\Delta(w)=1\otimes 1\otimes 1=({\rm id}\otimes\Delta)\Delta(w)$.

Assume that the result is true  for any $w$ with $({\rm dep}(w), {\rm bre}(w))< (p, q)$.

Let  $w\in \Phi_{\infty}(X)$ with $({\rm dep}(w), {\rm bre}(w))=(p,q)>(0,0)$. Then we have two cases to consider.

Case 1. If $q=1$, then $w=x\in X$ or $w=R(w')$ or $w=S(w')$.

Subcase  1.  If $w=x\in X$, then
$$
(\Delta\otimes {\rm id})\Delta(x)=x\otimes1\otimes1+1\otimes x\otimes1+1\otimes1\otimes x=({\rm id}\otimes\Delta)\Delta(x).
$$

Subcase  2. If  $w=S(w')$, then
\begin{eqnarray*}
&&({\rm id}\otimes\Delta)\Delta(S(w'))\\
&=&({\rm id}\otimes\Delta)(R(w')\otimes 1+({\rm id}\otimes S)\Delta (w'))\\
&=&  R(w')\otimes 1 \otimes 1+({\rm id}\otimes (\Delta S))\Delta (w'))\\
&=&  R(w')\otimes 1 \otimes 1+({\rm id}\otimes  R)\Delta (w')\otimes1+({\rm id}\otimes{\rm id}\otimes S)({\rm id}\otimes \Delta)\Delta(w')\\
&=&  R(w')\otimes 1 \otimes 1+({\rm id}\otimes  R)\Delta (w')\otimes1+({\rm id}\otimes{\rm id}\otimes S)(\Delta\otimes{\rm id})\Delta(w')\\
&=&  R(w')\otimes 1 \otimes 1+({\rm id}\otimes  R)\Delta (w')\otimes1+ (\Delta\otimes S)\Delta(w')\\
&=&  (R(w')\otimes 1  +({\rm id}\otimes  R)\Delta (w'))\otimes1+ (\Delta\otimes S)\Delta(w')\\
&=& \Delta(R(w'))\otimes 1+(\Delta\otimes S)\Delta(w')\\
&=&(\Delta\otimes{\rm id})(R(w')\otimes 1+({\rm id}\otimes  S)\Delta(w'))\\
&=&(\Delta\otimes{\rm id})\Delta (S(w')).
\end{eqnarray*}

Subcase  3. If $w=R(w')$,   then similar to the Subcase 2,   we have
$$
(\Delta\otimes {\rm id})\Delta(w)=({\rm id}\otimes\Delta)\Delta(w).
$$

Case 2. If $q>1$, then we can let $w=w_1w_2$ with $w_1, w_2\neq 1$.      Let us denote
$$
\Delta(w_1)=\sum_{(w_1)} w_{1(1)}\otimes w_{1(2)}, \ \ \
\Delta(w_2)=\sum_{(w_2)} w_{2 (1)}\otimes w_{2(2)}.
$$
Then
\begin{eqnarray*}
({\rm id}\otimes\Delta)\Delta(w)&=&({\rm id}\otimes\Delta)\Delta(w_1\diamond w_2)\\
&=& ({\rm id}\otimes\Delta)(\Delta(w_1)\diamond \Delta(w_2))\\
&=&\sum_{(w_1)}\sum_{(w_2)}  (w_{1(1)} \diamond w_{2(1)}) \otimes  \Delta (w_{1(2)}\diamond w_{2(2)})\\
&=&\sum_{(w_1)}\sum_{(w_2)}  (w_{1(1)} \diamond w_{2(1)}) \otimes  (\Delta (w_{1(2)})\diamond \Delta(w_{2(2)}))\\
&=&\sum_{(w_1)} \sum_{(w_2)} (w_{1(1)}\otimes \Delta (w_{1(2)})\diamond     (w_{2(1)}) \otimes  \Delta(w_{2(2)}). \\
\end{eqnarray*}
Similarly, we have
$$
(\Delta \otimes {\rm id})\bigtriangleup(w) =\sum_{(w_1)} \sum_{(w_2)} (\Delta(w_{1(1)})\otimes  w_{1(2)}) \diamond     (\Delta(w_{2(1)}) \otimes   w_{2(2)}).
$$
By  induction hypotheses, we have
$$
(\Delta\otimes {\rm id})\Delta(w_1)=({\rm id}\otimes\Delta)\Delta(w_1),\ \  (\Delta\otimes {\rm id})\Delta(w_2)=({\rm id}\otimes\Delta)\Delta(w_2).
$$
Therefore,
$$
(\Delta\otimes {\rm id})\Delta(w)=({\rm id}\otimes\Delta)\Delta(w).
$$
This completes the proof of coassociativity by induction.

We also  use   induction on $(m,n)$ for $w\in \Phi_{\infty}(X)$,    where $m={\rm dep}(w)$ and  $n={\rm bre}(w)$ to check the left  counicity conditions:
$$
(\varepsilon\otimes{\rm id})\Delta (w)=\beta_l(w),
$$
where
$$
\beta_l(w): RS(X)\rightarrow k\otimes RS(X), w\mapsto 1_k\otimes w.
$$

If $(m,n)=(0,0)$, then it is easy to see that the result is true.

Assume that the result is true for any $w$ with $(m, n)< (p, q)$.

Let  $w\in \Phi_{\infty}(X)$ with $({\rm dep}(w), {\rm bre}(w))=(p,q)>(0,0)$.

Case I. If ${q=\rm bre}(w)=1$, then $w=x\in X$ or $w=R(w')$ or $w=S(w')$.

Subcase I-1. If $w=x\in X$, then
$$
(\varepsilon\otimes{\rm id})\Delta (x)=1_k\otimes x=\beta_l(x).
$$

Subcase I-2. If  $w=S(w')$, then
\begin{eqnarray*}
(\varepsilon\otimes {\rm id})\Delta(w)&=&(\varepsilon\otimes {\rm id})\Delta (S(w'))\\
&=& (\varepsilon\otimes {\rm id})(R(w')\otimes 1+({\rm id}\otimes S)\Delta(w'))\\
&=& (\varepsilon\otimes {\rm id}) ({\rm id}\otimes S)\Delta(w'))\\
&=& ({\rm id}\otimes S)(\varepsilon\otimes {\rm id})\Delta(w'))\\
&=& ({\rm id}\otimes S)(1_k\otimes w')=1_k\otimes S(w')=\beta_l(w).
\end{eqnarray*}

Subcase I-3. If  $w=R(w')$, then similar to the Subcase I-2,  $(\varepsilon\otimes {\rm id})\Delta(R(w'))=\beta_l(R(w'))$.

Case II. If $q={\rm bre}(w)>1$, then we can let  $w=w_1w_2$ with $w_1, w_2\neq 1 $.     Let
$$
\Delta(w_1)=\sum_{(w_1)} w_{1(1)}\otimes w_{1(2)}, \ \ \
\Delta(w_2)=\sum_{(w_2)} w_{2 (1)}\otimes w_{2(2)}.
$$
By   induction  hypotheses, we have
$$
(\varepsilon\otimes{\rm id})\Delta (w_1)=\beta_l(w_1), \ \ (\varepsilon\otimes{\rm id})\Delta (w_2)=\beta_l(w_2).
$$
That is
$$
\sum_{(w_1)} \varepsilon(w_{1(1)})\otimes w_{1(2)}=1_k\otimes w_1, \ \ \sum_{(w_2)}\varepsilon( w_{ 2(1)})\otimes w_{2(2)}=1_k\otimes w_2.
$$
Thus
\begin{eqnarray*}
(\varepsilon\otimes {\rm id})\Delta(w)&=&(\varepsilon\otimes {\rm id})\Delta (w_1\diamond w_2)\\
&=& (\varepsilon\otimes {\rm id})(\Delta (w_1)\diamond \Delta(w_2))\\
&=& \sum_{(w_1), (w_2)} \varepsilon(w_{1(1)}\diamond w_{2 (1)}) \otimes (w_{1(2)} \diamond w_{2(2)})  \\
&=& \sum_{(w_1), (w_2)} (\varepsilon(w_{1(1)})\varepsilon(w_{ 2(1)})) \otimes (w_{1(2)} \diamond w_{2(2)})  \\
&=& \sum_{(w_1), (w_2)} (\varepsilon(w_{1(1)})\otimes  w_{1(2)})  (\varepsilon(w_{2 (1)})\diamond w_{2(2)}) \\
&=& (1_k\otimes w_1)(1_k\otimes w_2)= 1_k\otimes (w_1\diamond w_2)=1_k\otimes w=\beta_l(w).
\end{eqnarray*}
This completes the proof of $(\varepsilon\otimes{\rm id})\Delta (w)=\beta_l(w)$ by induction. \hfill $ \square$\\\\

{\bf Remark:}  Note that  $\varepsilon$ is not a right counit on $RS(X)$.  If $R\neq S$, then
$$
({\rm id}\otimes\varepsilon)\Delta(S(1))\neq \beta_r(S(1)),
$$
where
$$
\beta_r: RS(X)\rightarrow RS(X)\otimes k,  w\mapsto w \otimes 1_k.
$$\\

By Lemmas \ref{lem3.2}-\ref{lem3.4}, we have the following theorem.
\begin{theorem}
The quintuple $(RS(X), \diamond, u, \Delta, \varepsilon)$ is a left counital bialgebra.
\end{theorem}

\subsection{A left counital Hopf algebra structure}

For a $k$-algebra $(A, \mu, u)$ and a  left counital bialgebra coalgebra $(C, \Delta, \varepsilon)$, we define the convolution of two linear maps $f, g\in {\rm Hom}(C, A)$ to be the map
$
f\ast g\in {\rm Hom}( C, A)
$
given by the composition
$$
C\overset{\Delta}\longrightarrow  C\otimes C \overset{f\otimes g} \longrightarrow A\otimes A \overset{\mu}\longrightarrow  A.
$$
Let $(H, \mu, u, \Delta, \varepsilon)$ be a  left counital   bialgebra. A $k$-linear endomorphism $T$ of $H$ is called a right  antipode of  $H$ if it is the right inverse of ${\rm id}_H$ under the convolution product
$$
 {\rm id}_H\ast T=u\varepsilon.
$$
A left counital  Hopf algebra is a  left counital  bialgebra $H$ with an right antipode $T$. For more about  (left counital )  Hopf algebra, see for instance \cite{abe1980, zg2017, lguo}.

Recall that a     left counital bialgebra $(H, \mu, u, \Delta, \varepsilon)$ is called a graded  left counital bialgebra if there are a sequence of $k$-vector spaces $H^{(n)},\ n\geq 0,$  such that
\begin{enumerate}
\item[(a)]  $H=\bigoplus_{n=0}^{\infty}H^{(n)}$;
\item[(b)] For any $p, q\geq 0$,  $H^{(p)}H^{(q)}\subseteq H^{(p+q)}$;
\item[(c)] For any $n\geq 0$,  $\Delta(H^{(n)})\subseteq \bigoplus_{p+q=n}  H^{(p)}\otimes H^{(q)}$.
\end{enumerate}
A graded left counital  bialgebra  $H=\bigoplus_{n=0}^{\infty}H^{(n)}$ is called connected if  $H^{(0)}=k$  and $ker \varepsilon =\bigoplus_{n>0} H^{(n)}$.

Define $H_{RS}^{(n)}$ the $k$-linear space spanned by $\{w\in \Phi_{\infty}(X)| {\rm deg}(w)=n\}$, i.e.
$$
H_{RS}^{(n)}=k\{w\in \Phi_{\infty}(X)| {\rm deg}(w)=n\}.
$$
Then
$$
H_{RS}=\bigoplus_{n=0}^{\infty}H_{RS}^{(n)}, \ \ \ \  H_{RS}^{(0)}=k.
$$

\begin{lemma} \label{lem3.6}Let $p, q\geq 0$. Then
$$
H_{RS}^{(p)}\diamond H_{RS}^{(q)} \subseteq H_{RS}^{(p+q)}.
$$
\end{lemma}
{\bf Proof.} We just have  to prove that $w\diamond v\in H_{RS}^{(p+q)}$ for any  $w, v \in \Phi_{\infty}(X)$ with ${\rm deg}(w)=p$ and ${\rm deg}(v)=q$.  Induction on $n=p+q$.

If $n=0$, then $w=v=1$ and it is easy to see that the result is true.

Assume that the result is true for $n=p+q$.

Let  $w, v \in \Phi_{\infty}(X)$ with  ${\rm deg}(w)+{\rm deg}(v)=p+q+1$. Let
$$
w=w_1w_2\cdots w_t, v=v_1v_2\cdots v_m.
$$

Case 1. If $(w_t, v_1)\notin \bigcup _{Q\in \{R, \ S \}}Q(\Phi_{\infty}(X))\times Q(\Phi_{\infty}(X)) $, then $w\diamond v=wv\in \Phi_{\infty}(X)$ and it is easy to see the result is true.

Case 2.  If $(w_t, v_1)\in \bigcup _{Q\in \{R, \ S \}}Q(\Phi_{\infty}(X))\times Q(\Phi_{\infty}(X)) $, then $w_t=R(w'), v_1=R(v')$ or $w_t=S(w'), v_1=S(v')$.

(i) If $w_t=S(w'), v_1=S(v')$, then
\begin{eqnarray*}
w\diamond v&=&w_1w_2\cdots w_{t-1}(S(w')\diamond S(v'))v_2\cdots v_m\\
&=& w_1w_2\cdots w_{t-1}S(R(w')\diamond v'+w'\diamond S(v'))v_2\cdots v_m.
\end{eqnarray*}

Let $l={\rm deg}(w_t)+{\rm deg}(v_1)$.
By  induction hypotheses,  $R(w')\diamond v', w'\diamond S(v')\in H_{RS}^{(l-1)}$. Therefore, $w\diamond v\in H_{RS}^{(p+q+1)}$.

(ii) If $w_t=R(w'), v_1=R(v')$, then similar to the proof of   (i), we have the result is true.

This completes the proof by induction. \hfill $ \square$\\

\begin{lemma} \label{lem3.7}For any $n\geq 0$,
$$
\Delta(H_{RS}^{(n)})\subseteq \bigoplus_{p+q=n}H_{RS}^{(p)}\otimes H_{RS}^{(q)}.
$$
\end{lemma}
{\bf Proof.}
For  $n=0, 1$, it is easy to see that the result is true.   Assume that the result is true for $0 \leq n\leq m$.
We just have to prove that $\Delta(w)\in \bigoplus_{p+q=m+1}H_{RS}^{(p)}\otimes H_{RS}^{(q)}$ for any $w$ with ${\rm deg}(w)=m+1$.

If ${\rm bre}(w)=1$, then we have $w=R(w')$ or $w=S(w')$.

If $w=S(w')$, then by induction hypotheses, we have \
$$
\Delta(w)=\Delta(S(w'))=R(w')\otimes 1+({\rm id}\otimes S)\Delta(w')\in \bigoplus_{p+q=m+1}H_{RS}^{(p)}\otimes H_{RS}^{(q)}.
$$

Similarly,  if $w=R(w')$, then $\Delta(w) \in \bigoplus_{p+q=m+1}H_{RS}^{(p)}\otimes H_{RS}^{(q)}$.

If ${\rm bre}(w)\geq 2$,  say $w=uv$ with $u, v\neq 1$, then by induction hypotheses and Lemma \ref{lem3.6}, we have
\begin{eqnarray*}
\Delta(w)&=&\Delta(u)\diamond \Delta(v)\\
&\in& (\bigoplus_{p+q={\rm deg}(u)}H_{RS}^{(p)}\otimes H_{RS}^{(q)} )\diamond(\bigoplus_{p'+q'={\rm deg}(v)}H_{RS}^{(p')}\otimes H_{RS}^{(q')} )\\
&\subseteq& \bigoplus_{t+s=m+1}H_{RS}^{(t)}\otimes H_{RS}^{(s)}.
\end{eqnarray*}
 This completes the proof by induction. \hfill $ \square$ \\

By a similar proof of the theorems in  (\cite{zg2017,lguo}), we have

\begin{lemma}\label{th3.8}
A connected left counital  bialgebra $(H, \mu, u, \Delta, \varepsilon)$ is a  left counital  Hopf algebra.
\end{lemma}

By Lemmas \ref{lem3.6}-\ref{th3.8}, we have the following theorem.
\begin{theorem}
The free Rota-Baxter system  $RS(X)=\bigoplus_{n=0}^{\infty}H_{RS}^{(n)}$ is  a connected  left counital  bialgebra. It follows that  $RS(X)$ is a left counital  Hopf algebra.
\end{theorem}

\ \

\noindent{\bf Acknowledgement}
 We wish to express our thanks to Dr. Xing Gao and Prof. Li Guo for helpful and valuable suggestions and comments.

\end{document}